\documentclass[11pt,a4paper]{article}
\usepackage{graphicx}
\usepackage{amsfonts}
\usepackage{latexsym}
\newcommand{\beql}[1]{\begin{equation}\label{#1}}
\newcommand{\eeq}{\end{equation}}
\newcommand{\comment}[1]{}

\newcommand{\eqref}[1]{{\rm (\ref{#1})}}

\newcommand{\Abs}[1]{{\left|{#1}\right|}}

\newcommand{\Qed}{\ \\\mbox{$\Box$}}

\newcommand{\Set}[1]{{\left\{{#1}\right\}}}

\newcommand{\RR}{{\Bbb R}}

\newcommand{\ZZ}{{\Bbb Z}}

\newcommand{\one}{{\bf 1}}

\newcommand{\inner}[2]{{\langle #1, #2 \rangle}}

\newcommand{\dens}{{\rm dens\,}}

\newcommand{\ft}[1]{\widehat{#1}}

\newcounter{open}
\newcounter{dfn}
\def\thedfn{\arabic{dfn}}

\newcounter{obs}
\def\theobs{\arabic{obs}}

\newcounter{thm}
\newcounter{mysec}
\def\themysec{\arabic{mysec}}
\newcommand{\mysection}[1]{
  \vskip 0.25in
  \refstepcounter{mysec}\noindent{\large\bf \S\themysec.\ {#1}}\par
  \nopagebreak
  \addcontentsline{toc}{section}{{\bf \themysec.}\ {#1}}
}

\newcounter{mysubsec}[mysec]
\newtheorem{theorem}{Theorem}

\begin{document}

\begin{center}
{\Large \bf The Steinhaus tiling problem and the range of certain quadratic forms}\\
\ \\
{\sc Mihail N. Kolountzakis} and {\sc Michael Papadimitrakis}\\
\ \\
\small September 2000
\end{center}

\centerline{\em Dedicated to the memory of Tom Wolff}

\begin{abstract}
We give a short proof of the fact that there are no measurable
subsets of Euclidean space (in dimension $d\ge 3$), which, no matter
how translated and rotated, always contain exactly one integer lattice
point.
In dimension $d=2$ (the original Steinhaus problem) the question
remains open.
\end{abstract}

\mysection{The Steinhaus problem}

Steinhaus (1957) asked if there exists a subset of the plane which,
no matter how translated and rotated, always contains exactly one
point with integer coordinates.
This question remains unanswered.

In this paper we deal only with the measurable version of the Steinhaus
problem in dimension $d\ge2$,
in which a {\em measurable} subset $E$ of $\RR^d$ is sought
with the property that for almost every $x \in \RR^d$ and for almost every
isometry ${\cal S}:\RR^d\to\RR^d$
$$
\Abs{({\cal S}E+x) \cap \ZZ^d} = 1.
$$
The {\em Steinhaus property} may also be written as follows.
For almost all isometries ${\cal S}$
$$
\sum_{n\in\ZZ^d} \one_{{\cal S}E}(x - n) = 1, \ \ \mbox{a.e.($x$)}.
$$
In other words, the set $E$ must be such that almost all of its rotations
${\cal S}E$ must {\em tile} when translated at the locations $\ZZ^d$.

In a recent paper \cite{KW} Wolff showed that there are no
Steinhaus sets in dimension $d\ge 3$.
He showed much more:
if $f \in L^1(\RR^d)$, $d\ge3$, is a {\em Steinhaus function}, i.e., if
$$
\sum_{n\in\ZZ^d} f(x - {\cal S}n) = 1, \ \ \mbox{a.e.($x$)},
$$
for a dense set of isometries ${\cal S}$, then $f$ is almost everywhere
equal to a continuous function.
This, of course, implies that no Steinhaus {\em sets} exist for $d\ge 3$.
We refer the reader to \cite{KW} for additional results regarding properties
of Steinhaus sets in dimension $2$ (if they exist) and references to other work.

Suppose that $\Lambda = A \ZZ^d \subset \RR^d$, $A \in GL(d, \RR)$,
is a lattice and let $\Lambda^* = A^{-\top}\ZZ^d$ be its dual lattice.
By elementary harmonic analysis one can see that,
for an $L^1$ function $f$, we have
$$
\sum_{\lambda\in\Lambda} f(x-\lambda) = C, \ \ \mbox{a.e.($x$)},
$$
if and only if its Fourier Transform $\ft{f}$ vanishes on
$\Lambda^*\setminus\Set{0}$.
Integrating over a large region
it is easy to see that the constant $C$ is equal to the density
of the lattice $\Lambda$ times the integral of $f$.

It follows that for $E$ to be a Steinhaus set it is necessary and
sufficient that $\Abs{E}=1$ and that
$\ft{\one_E}$ vanishes on all rotations of
the (self-dual) lattice $\ZZ^d$, except at $0$,
i.e., it is necessary and sufficient
that $\ft{\one_E}$ vanishes on all spheres
with positive radius centered at the origin
which go through at least one integer lattice point.

In this paper we will show that there are no Steinhaus sets
in dimension $d\ge 3$.
The method relies on some arithmetic properties of certain
quadratic forms in $d$ variables and is overall much simpler
than the method used in \cite{KW}.
There, of course, much stronger results were proved, using
advanced methods of harmonic analysis, about Steinhaus {\em functions}.
As mentioned above, these results have as a corollary the non-existence
of Steinhaus sets for $d\ge 3$.
Our method does not seem capable of giving any interesting results
about Steinhaus functions.
(These do exist: take any $L^1$ function whose Fourier Transform vanishes
on all spheres centered at the origin that go through a lattice
point.)

Our method is not applicable for the $d=2$, and we include a proof
of this.

The case $d\ge 4$ is presented separately from $d=3$ (from which it
follows) since it is much simpler.

\noindent
{\bf Acknowledgment.}
We are indebted to Professors A. Bremner and N. Tzanakis for very
valuable suggestions.

\mysection{The key observation}

In any dimension $d$ write ${\cal B}$ for the union of all spheres
centered at the origin that go through at least one lattice point.
The point $0$ is included in ${\cal B}$.

Assume from now on that the set $E$ is a Steinhaus set in dimension $d$.

Suppose now that we can find a lattice $\Lambda^* \subset {\cal B}$
with $\det\Lambda^*$ not an integer.
Since $\ft{\one_E}$ vanishes on $\Lambda^*\setminus\Set{0}$ it
follows that $E + \Lambda$ is a tiling at level
$C = \Abs{E}\times\dens\Lambda = 1 \times \det\Lambda^*$,
which is not an integer.
This is a contradiction as, obviously, any set may only
tile at an integral level.

Hence, there are no Steinhaus sets in dimension $d$ if one can
find a lattice of non-integral volume which is contained in
${\cal B}$.
Since a point $x\in\RR^d$ belongs to ${\cal B}$ if and only if
$\Abs{x}^2$ is a sum of $d$ integer squares, we obtain the
following Theorem, by looking at the quadratic form $\inner{A^\top Ax}{x}$
for each lattice $\Lambda^* = A\ZZ^d$.
\begin{theorem}
\label{th:main}
If there exists a positive definite quadratic form
$Q(x) = Q(x_1,\ldots,x_d) = \inner{Bx}{x}$ such that for all integral
$x_1,\ldots,x_d$ its value is the sum of $d$ integer squares,
and the determinant of $Q$, $\det B$, is not the square of an integer,
then there are no Steinhaus sets in dimension $d$.
\end{theorem}

\mysection{Dimension $d \ge 4$}

Consider the $4\times 4$ matrix $B$ with $1$ on the diagonal
and $1/2$ everywhere else.
The matrix $B$ is positive definite (its eigenvalues are
$1/2$, $1/2$, $1/2$ and $5/2$) and its determinant is
$5/16$.
It defines the quadratic form
$$
Q(x) = Q(x_1, \ldots, x_4) = \inner{Bx}{x} = \sum_{i=1}^4 x_i^2 +
	\sum_{i>j} x_i x_j,
$$
which is obviously integer valued and has non-square determinant.
Furthermore, every non-negative integer may be written as a sum
of four squares (Lagrange).
It follows from Theorem \ref{th:main} that there are no
Steinhaus sets for $d=4$.

The same is true for dimension $d>4$ as one may consider the matrix
which has $B$ in its upper left $4\times 4$ corner and is equal to the
identity matrix elsewhere.

\mysection{Dimension $d=3$}

The determinant of the form that appears in the following Theorem
is $2 \cdot 11 \cdot 6$, which is not a square,
hence there are no Steinhaus sets in dimension $3$.
\begin{theorem}
\label{th:3}
For each $x, y, z \in \ZZ$ the number
$$
2x^2 + 11y^2 + 6z^2
$$
is a sum of three integer squares.
\end{theorem}
{\bf Proof.}
Suppose this is false and that there are $(x_0, y_0, z_0) \neq (0, 0, 0)$
and
\begin{itemize}
\item[(a)]
$Q(x_0, y_0, z_0)$ is not a sum of three squares, and
\item[(b)]
$x_0^2 + y_0^2 + z_0^2$ is minimal.
\end{itemize}
From (a), and the well known characterization of those natural numbers
that cannot be written as a sum of three squares, we have that
$$
Q(x_0, y_0, z_0) = 4^\nu(8k+7),\ \ \nu\ge0, k\ge0.
$$
If all $x_0, y_0, z_0$ are even, we have $\nu\ge 1$, and, setting
$x_0 = 2 x_1$, $y_0 = 2y_1$ and $z_0 = 2z_1$,
we obtain that $Q(x_1, y_1, z_1)$ is not a sum of three squares, which
contradicts the minimality of the initial triple $(x_0, y_0, z_0)$.
We conclude that at least one of $x_0, y_0, z_0$ is odd.

\noindent\underline{\bf Case No 1:} $\nu=0$.

Then $Q(x_0, y_0, z_0) = 7 \bmod 8$.
But the quadratic residues mod $8$ are $0$, $1$ and $4$, and one
checks by examining all the possibilities that $Q$ is never $7 \bmod 8$.

\noindent\underline{\bf Case No 2:} $\nu=1$.

Then $Q(x_0, y_0, z_0) = 32k + 28$.
Hence $y_0$ is even, say $y_0 = 2 y_1$.
We get
$$
x_0^2 + 22y_1^2 + 3 z_0^2 = 16k + 14,
$$
from which we conclude that $x_0$ and $z_0$ are odd,
$x_0 = 2x_1 + 1$, $z_0 = 2z_1 + 1$.
Substitution gives
\begin{eqnarray*}
4x_1^2 + 4x_1 + 1 + 22y_1^2 + 12 z_1^2 + 12 z_1 + 3 &=& 16k + 14\\
2x_1(x_1+1) + 11y_1^2 + 6z_1(z_1+1) + 2 &=& 8k+7\\
2x_1(x_1+1) + 11y_1^2 + 6z_1(z_1+1) &=& 5 \bmod 8.
\end{eqnarray*}
But $\xi^2+\xi = 0$ or $2$ or $4$ or $6 \bmod 8$, for all $\xi$,
hence, by applying this to the first and last term in the above sum,
and checking all possibilities we get a contradiction.

\noindent\underline{\bf Case No 3:} $\nu\ge 2$.

As in Case No 2: $y_0 = 2y_1$, $z_0 = 2z_1+1$, $x_0 = 2x_1+1$.
Hence
$$
2x_1(x_1+1) + 11y_1^2 + 6z_1(z_1+1) + 2 = 4^{\nu-1}(8k+7),\ \ \nu-1\ge1.
$$
So $y_1$ is even, $y_1 = 2 y_2$, which gives
$$
x_1(x_1+1) + 22y_2^2 + 3 z_1(z_1+1) + 1 = 2\cdot 4^{\nu-2}(8k+7),
$$
a contradiction as the left hand side is odd while the right hand side
is even.
\Qed

\mysection{Dimension $d=2$}

Our method cannot give any results in dimension $2$:
\begin{theorem}
\label{th:2}
Any positive-definite binary quadratic form whose values are always sums of two
integer squares must have a determinant which is the square of an integer.
\end{theorem}
{\bf Proof.}
Let $Q(x, y)$ be such a quadratic form, which we may write
as
$$
Q(x, y) = \Abs{A\left(\begin{array}{c}x \\ y\end{array}\right)}^2,
$$
where $A = \left(\begin{array}{cc}a & c \\ b & d\end{array}\right)$.
In this notation the determinant of $Q$ is $(\det A)^2$.
We now use the following theorem.
\begin{theorem}
\label{th:dls}
{\rm(Davenport, Lewis and Schinzel \cite{DLS})}
Suppose that $f \in \ZZ[t]$ is such that for each arithmetic progression
$S$ there is $t\in S$ such that $f(t)$ is a sum of two squares.
Then there are $x, y \in \ZZ[t]$, such that
$$
f(t) = x^2(t) + y^2(t).
$$
\end{theorem}
{\bf Remarks.}

1. The assumptions of Theorem \ref{th:dls} are much weaker than we can afford.

2. Theorem \ref{th:dls} was proved to answer a question raised by LeVeque,
who had asked if any polynomial in $\ZZ[t]$ whose values are
always sums of two squares must be a sum of squares of two
linear forms with integer coefficients.

Let $f(t) = Q(t, 1)$, which is a polynomial with integer coefficients.
It follows that there are integers $\alpha, \beta, \gamma, \delta$
such that, for all $t\in\ZZ$,
$$
(at+c)^2 + (bt+d)^2 = (\alpha t + \gamma)^2 + (\beta t + \delta)^2.
$$
Expanding, and identifying the coefficients we obtain
$$
a^2+b^2 = \alpha^2+\beta^2,\ ac+bd = \alpha\gamma+\beta\delta,\ \mbox{and }
	c^2+d^2 = \gamma^2+\delta^2.
$$
We have for the determinant of $Q$:
\begin{eqnarray*}
\det Q &=& (ad-bc)^2\\
	&=& \Abs{\begin{array}{cc}a^2+b^2 & ac+bd \\ ac+bd & c^2+d^2\end{array}}\\
	&=& \Abs{\begin{array}{cc}\alpha^2+\beta^2 & \alpha\gamma+\beta\delta \\
		\alpha\gamma+\beta\delta & \gamma^2+\delta^2\end{array}}\\
	&=& (\alpha\delta-\beta\gamma)^2,
\end{eqnarray*}
which is the square of an integer.
\Qed

\noindent
\ \\

\noindent
{\sc\small
Department of Mathematics, University of Crete, Knossos Ave.,
714 09 Iraklio, Greece.\\
E-mail: {\tt mk@fourier.math.uoc.gr}, {\tt papadim@math.uoc.gr}
}

\end{document}